\documentclass[12pt,reqno,pdftex]{amsart}
\usepackage{amsmath}
\usepackage{amsthm}
\usepackage{graphicx}
\usepackage{bbm,amssymb}
\usepackage[hyperfootnotes=false,colorlinks]{hyperref}

\newcommand{\ceiling}[1]{\left\lceil {#1} \right\rceil}

\newtheorem{theorem}{Theorem}[section]
\newtheorem{prop}[theorem]{Proposition}
\newtheorem{lemma}[theorem]{Lemma}

\newtheorem{conjecture}[theorem]{Conjecture}

\theoremstyle{remark}
\newtheorem*{remark}{Remark}

\theoremstyle{definition}
\newtheorem*{definition}{Definition}

\def\xx{\mathbf{x}}
 
\DeclareSymbolFont{AMSb}{U}{msb}{m}{n}
\DeclareMathSymbol{\C}{\mathbin}{AMSb}{"43} 
\DeclareMathSymbol{\EE}{\mathbin}{AMSb}{"45} 
\DeclareMathSymbol{\N}{\mathbin}{AMSb}{"4E} 
\DeclareMathSymbol{\PP}{\mathbin}{AMSb}{"50} 
\DeclareMathSymbol{\Q}{\mathbin}{AMSb}{"51} 
\DeclareMathSymbol{\R}{\mathbin}{AMSb}{"52} 
\DeclareMathSymbol{\Z}{\mathbin}{AMSb}{"5A}

\begin{document}
%\title{Sandpiles: Growth Rates and Explosions}
%\title{Diameter Bounds and Explosions in Sandpiles}
%\title{Sandpile Growth Rates and Explosions}
\title{Growth Rates and Explosions in Sandpiles}
%\title[Growth Rate and Stabilizability of Sandpiles]{On the Growth Rate and Stabilizability of Sandpiles in $\Z^d$}
\author{Anne Fey \and Lionel Levine \and Yuval Peres}

\address{Anne Fey, Delft Institute of Applied Mathematics, Delft University of Technology, The Netherlands, {\tt a.c.fey-denboer@tudelft.nl}} \address{Lionel Levine, Department of Mathematics, Massachusetts Institute of Technology, Cambridge, MA 02139, {\tt \url{http://math.mit.edu/~levine}}}  \address{Yuval Peres, Theory Group, Microsoft Research, Redmond, WA 98052, {\tt peres@microsoft.com}}

\date{September 9, 2009}

\begin{abstract}
We study the abelian sandpile growth model, where~$n$ particles are added at the origin on a stable background configuration in $\Z^d$.  Any site with at least $2d$ particles then topples by sending one particle to each neighbor.
We find that with constant background height $h \leq 2d-2$, the diameter of the set of sites that topple has order~$n^{1/d}$. This was previously known only for $h < d$.  Our proof uses a strong form of the least action principle for sandpiles, and 
a novel method of background modification.  

We can extend this diameter bound to certain backgrounds in which an arbitrarily high fraction of sites have height $2d-1$.  On the other hand, we show that if the background height $2d-2$ is augmented by $1$ at an arbitrarily small fraction of sites chosen independently at random, then adding finitely many particles creates an \emph{explosion} (a sandpile that never stabilizes).
\end{abstract}

\maketitle

\section{Introduction}

\begin{figure}
\centering
\includegraphics[width=.48\textwidth]{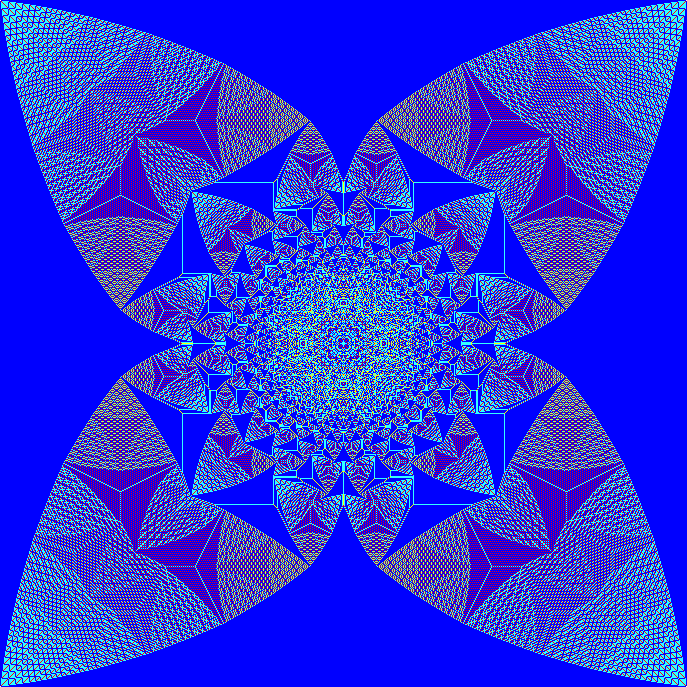} \hspace{0.7mm}
\includegraphics[width=.48\textwidth]{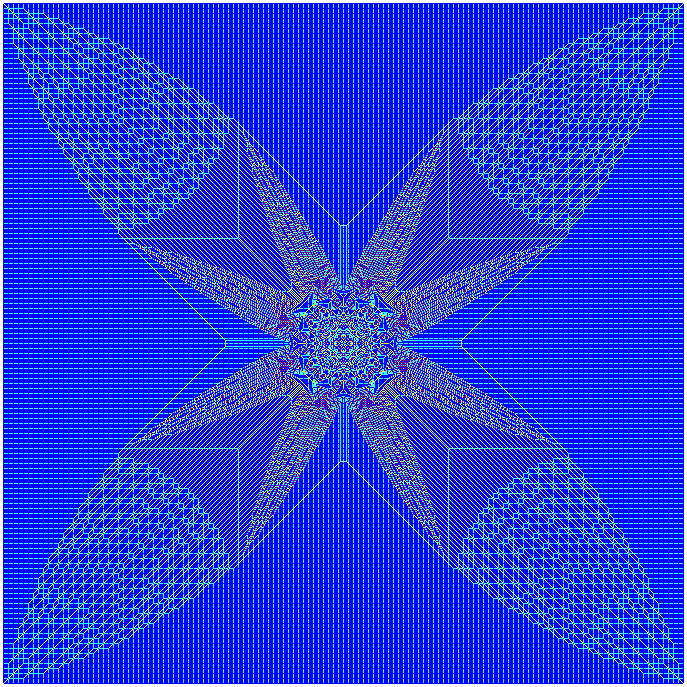}
\caption{Left: Stable sandpile of $n=2 \cdot 10^5$ particles in~$\Z^2$ on background height $h=2$.  Right: Sandpile of $n=15000$ particles in $\Z^2$ on background height~$3$, except every fifth row and column has background height~$2$.  In both cases, the set $T_n$ is a square.
Color scheme: sites colored blue have~$3$ particles, turquoise~$2$ particles, yellow~$1$ particle, red~$0$ particles.}
\label{fig:boxes}
\end{figure}

In this paper we consider the abelian sandpile model as a growth model in the integer lattice $\Z^d$.
The model starts from a stable {\em background} configuration in which each site $x$ has a pile of $\sigma(x) \leq 2d-1$ particles.  To this background, $n$ particles are added at the origin. Typically, $n$ is large. We {\em stabilize} this configuration by {\em toppling} every unstable site; that is, every site with at least $2d$ particles gives one particle to each of its neighbors, until there are no more unstable sites.
For more information on the abelian sandpile model, also known as the chip-firing game, see \cite{BTW,BLS,dhar,MRZ}.

To keep things simple in this introduction, we will enumerate the sites in $\Z^d$ (say, in order of increasing distance from the origin, breaking ties arbitrarily) and perform topplings one by one in discrete time: At each time step, if there are any unstable sites, then the smallest unstable site topples.  All of our results hold also for the more general toppling procedures discussed in section~\ref{sec:leastaction}.

Let $T_n = T_{n,d,\sigma}$ be the set of sites that topple (Figure~\ref{fig:boxes}).  Since these sets are nested, $T_1 \subseteq T_2 \subseteq \ldots$, it is natural to view them as a growth model, with~$n$ playing the role of a time parameter.  We distinguish between two extreme cases.  If~$T_n$ is finite for all~$n$, we say that $\sigma$ is \emph{robust}.  In this case we are interested in the growth rate, i.e.\ in how the diameter of~$T_n$ grows with~$n$.
%, as well as in the geometry of the set $T_n$.

At the other extreme, if $T_n = \Z^d$ for some $n$, then every site topples infinitely often.   Otherwise, some site $x \in \Z^d$ must finish toppling before all of its neighbors do; since each neighbor topples at least once after~$x$ finishes toppling,~$x$ receives~$2d$ additional particles and must topple again.  

If $T_n = \Z^d$ for some~$n$, then we say that $\sigma$ is \emph{explosive}, and it is \emph{exploding} when the~$n$ particles are added. (In \cite{feyredig}, the term `not stabilizable' was used for `exploding,'  and `metastable' for `explosive.') The simplest example of an explosive background is $\sigma(x)=2d-1$ for all~$x$, to which the addition of a single extra particle causes every site in $\Z^d$ to topple infinitely many times.  
%An exploding sandpile started from a different background is shown in Figure~\ref{fig:latticebootstrap}.

\begin{figure}
\centering
\includegraphics[scale=.25]{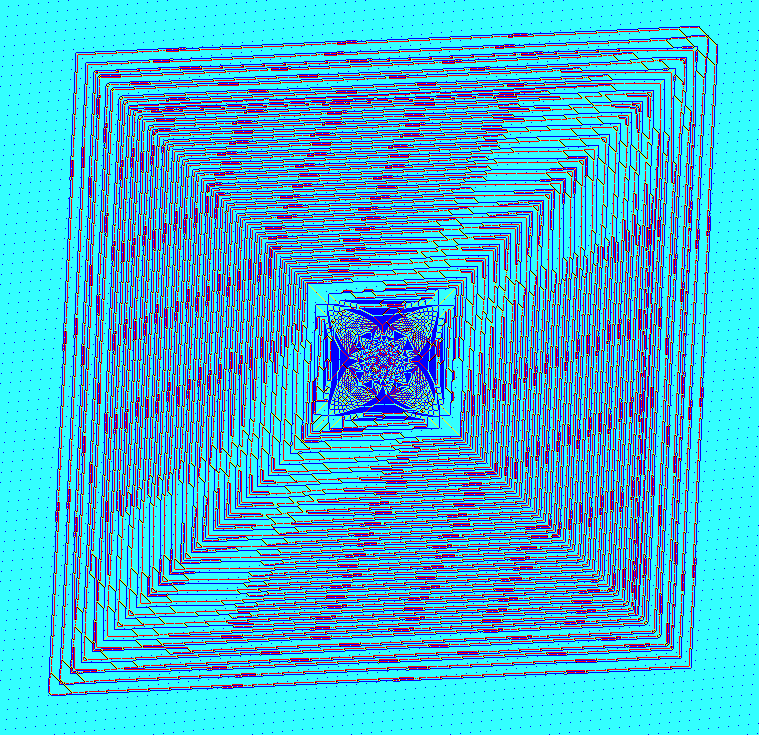}
\caption{An exploding sandpile started from $n=5000$ particles in $\Z^2$.  Background height is $2$ except for sites in the lattice generated by $(1,10)$ and $(10,1)$, which have background height $3$.  Unstable sites are colored black.}
\label{fig:latticebootstrap}
\end{figure}

We remark that an intermediate behavior is possible, when $T_n$ is infinite for a finite $n$, but $T_n \neq \Z^d$ for all $n$.  An example is the background with $3$ particles at every site on the $x$-axis in $\Z^2$, and $2$ particles at every other site.  Adding one particle at the origin produces an infinite avalanche of topplings, but each site topples only finitely many times.  This example shows that exploding is a strictly stronger condition than having an infinite avalanche of topplings.

%In this case we say that $\sigma$ is \emph{$\omega$-stable}. 
%As an example, consider the background on $\Z$ where every site $x \leq 0$ has one particle, and all other sites are empty. Then $T_n = (-\infty,n-1]$.
%every addition to the origin causes all sites to the left of the origin to topple, but only finitely many to the right of the origin.

The papers~\cite{shapes} and~\cite{LP09a} investigated 
the case of a robust constant background of $h \leq 2d-2$ particles at every site. In the regime $h<d$, the diameter of $T_n$ grows like $n^{1/d}$; the best known bounds can be found in \cite[Theorem~4.1]{LP09a}.  
%In a related paper \cite{dharpattern}, the focus is on the intriguing height patterns that are present in the stabilized configuration.
% not sure if this reference belongs here - they mainly studied other lattices
In the case $h=2d-2$, the set $T_n$ is a cube for every~$n$, and an upper bound for the radius is $n$ \cite[Theorem~4.1]{shapes}. No proof was found for a better upper bound, even though simulations clearly indicated a growth rate proportional to $n^{1/d}$.  

In this paper we complete the picture by deriving an upper bound of order~$n^{1/d}$ on the diameter of $T_n$ for all $h \leq 2d-2$, and even for some backgrounds arbitrarily close to $2d-1$.  We first correct a gap in the proof of the outer bound of \cite[Theorem~4.1]{LP09a}, which we thank Haiyan Liu for pointing out to us.  Then we use this theorem together with a new technique of ``background modification'' to extend the bounds to higher values of~$h$.

Throughout the paper, we will typically use the symbol $\sigma$ to indicate a stable background configuration, and $\eta$ to indicate an arbitrary (possibly unstable) configuration.  We use~$\overline{h}$ to denote the constant configuration $\sigma(x) \equiv h$, and we denote a single particle at the origin by~$\delta_o$.  Write
    \[ Q_r = \{x\in \Z^d \,:\, \max |x_i| \leq r \} \]
for the cube of side length $2r+1$ centered at the origin in $\Z^d$.
Let $\omega_d$ be the volume of the unit ball in $\R^d$.

%In this introductory section, for simplicity we state our results in dimension two.
Our main result, proved in section~\ref{growthratessection}, is the following.

%\begin{theorem}
%\label{thecubeintro}
%Let $T_{n,2,\overline{2}}$ be the set of sites in $\Z^2$ that topple during stabilization, after $n$ particles are added at the origin to a background of $2$ particles at every site in $\Z^2$.  Then for any $\epsilon>0$, we have
%	\[ T_{n,2,\overline{2}} \subset Q_r \]
%for all sufficiently large $n$, where
%	\[ r =  \left(2 + \epsilon \right) \sqrt{\frac{n}{\pi}}. \]
%\end{theorem}

\begin{theorem}
\label{thecube}
Fix integers $d \leq h \leq 2d-2$, and let $T_{n,d,\overline{h}}$ be the set of sites in $\Z^d$ that topple during the stabilization of $\overline{h}+n\delta_o$.
% if $n+h$ particles start at the origin and~$h$ particles start at every other site in $\Z^d$.  
Then for any $\epsilon>0$, we have
	\[ T_{n,d,\overline{h}} \subset Q_r \]
for all sufficiently large $n$, where
	\[ r =  \frac{d+\epsilon}{2d-1-h} \left( \frac{n}{\omega_d} \right)^{1/d}. \]
\end{theorem}

In the case $d=h=2$, Theorem~\ref{thecube} gives a bound of $\frac{2+\epsilon}{\sqrt{\pi}} \sqrt{n} \approx 1.13 \sqrt{n}$ on the radius of the square of sites that topple.  Large scale simulations by David Wilson indicate that the actual radius is approximately $0.75 \sqrt{n}$.

It is natural to ask what happens when the background height $h$ exceeds $2d-2$.  While the background $\overline{2d-1}$ is explosive, our next result shows that there exist robust backgrounds in which an arbitrarily high proportion of sites have $2d-1$ particles.  For $m\geq 1$, let
	\[ \Lambda(m) = \{x \in \Z^d : m \not|~x_i \mbox{ for all } 1\leq i\leq d \}. \]
Thus $\Lambda(m)$ is a union of cubes of side length $m-1$.   The following theorem generalizes the case $h=2d-2$ of Theorem~\ref{thecube}, which corresponds to the case $m=1$.

\begin{theorem}
\label{rowsandcolumns}
For any $m \geq 1$, the background
	\[ \sigma = \overline{2d - 2} + 1_{\Lambda(m)} \]
is robust on $\Z^d$.  Moreover, writing $T_{n,d,\sigma}$ for the set of sites in $\Z^d$ that topple during the stabilization of $\sigma+n\delta_o$, then for any $\epsilon>0$, we have
	\[ T_{n,d,\sigma} \subset Q_r \]
for all sufficiently large $n$, where
	\[ r =  m (d + \epsilon) \left( \frac{n}{\omega_d} \right)^{1/d}. \]
\end{theorem}

On the basis of this theorem, one might guess that $2d-1$ is the critical density below which a background is robust and above which it is explosive.  Our next two results show that this is not the case.  Starting from background height $2d-2$,  we can destroy robustness by adding extra particles on an arbitrarily sparse lattice $L \subset \Z^d$ (Figure~\ref{fig:latticebootstrap}).
 
 \begin{prop}
\label{latticebootstrapintro}
Let $\xx_i = (x_{i1}, \ldots, x_{id})$ for $i=1,\dots,d$ be linearly independent vectors in $\Z^d$ satisfying $\gcd(x_{1j},\ldots,x_{dj})=1$ for all $j=1,\ldots,d$.  Let $L = \Z \xx_1 + \ldots + \Z \xx_d$.  Then the background $\overline{2d-2}+1_L$ on $\Z^d$ is explosive.
\end{prop}

Comparing Theorem~\ref{rowsandcolumns} and Proposition~\ref{latticebootstrapintro}, we see that the geometry of the extra particles plays a more important role in determining robustness than the density of sites at which particles are added.  In particular, the fact that $L$ intersects every hyperplane in $\Z^d$ that is parallel to one of the coordinate hyperplanes, while $\Lambda(m)$ does not, plays a key role in the proofs.

As our next result shows, the lattice structure is not essential in Proposition~\ref{latticebootstrapintro}.  We can also produce an explosive background by adding particles at rare random sites.

\begin{prop}
\label{randombootstrapintro}
Fix $\epsilon>0$, and let $(\beta(x))_{x \in \Z^d}$ be independent Bernoulli random variables with $\PP(\beta(x)=1)=\epsilon$.  With probability $1$, the background $\overline{2d-2}+\beta$ on $\Z^d$ is explosive.
\end{prop}

Our proofs make extensive use of the abelian property of the abelian sandpile model, which we state and generalize in the next section.

\section{Least Action Principle}
\label{sec:leastaction}

We begin by recalling the notion of toppling procedure defined in~\cite{FMR}.  This formalism includes most of the natural ways to topple, including: discrete time parallel updates, in which all unstable vertices topple simultaneously; toppling in nested volumes, in which we successively stabilize larger and larger finite regions of $\Z^d$; and Markov toppling in continuous time, in which each site has a Poisson clock and attempts to topple whenever its clock rings.  The technical details of the toppling procedures are tangential to our main argument, so the reader may wish to skim them and move on to the ``\hyperref[odometerisminimal]{least action principle},'' which is the only new material in this section.

Let $\mathcal{X} = \Z^{\Z^d}$.  We think of elements of $\mathcal{X}$ as particle configurations on $\Z^d$ in which some sites may have a negative number of particles.  We endow $\mathcal{X}$ with the Borel $\sigma$-algebra coming from the product topology, with $\Z$ having the discrete topology.  On $\Z^d$ and on $\N$ we use the full power set as a $\sigma$-algebra, and on the half-line $[0,\infty)$ we use the usual Borel $\sigma$-algebra.  A \emph{toppling procedure} is a measurable function
%define \sigma-fields
	\[ T : [0,\infty) \times \Z^d \times \mathcal{X} \to \N \]
satisfying for all $\eta \in \mathcal{X}$ and all $x \in \Z^d$
\begin{itemize}
\item[(a)] $T(0,x,\eta) = 0$.
\item[(b)] The function $t \to T(t,x,\eta)$ is right-continuous and nondecreasing with jumps of size at most one, i.e., for all $t \geq 0$,
	\[ T(t,x,\eta) - \lim_{s\uparrow t} T(s,x,\eta) \leq 1. \]
\item[(c)] In every finite time interval, there are only finitely many jumps at~$x$.
\item[(d)] There is no ``infinite backward chain of topplings,'' i.e., no path $x_1 \sim x_2 \sim \ldots$ and sequence of times $t_1 > t_2 > \ldots$ such that for all $i=1,2,\ldots$
	\[ T(t_i,x_i,\eta) > \lim_{s\uparrow t_i} T(s,x_i,\eta). \]
\end{itemize}

We interpret $T(t,x,\eta)$ as the number of times~$x$ topples in the time interval~$[0,t]$, for initial configuration~$\eta$.  We say that~$x$ topples at time~$t$ for initial configuration $\eta$, if 	
	\[ T(t,x,\eta)> \lim_{s\uparrow t} T(s,x,\eta). \]
%Note that condition (d) holds automatically for discrete time toppling procedures.
The toppling procedure $T$ may be deterministic (as in parallel updates or nested volumes) or random (as in Markov toppling).  A random toppling procedure can be viewed as a measurable function
	\[ T : [0,\infty) \times \Z^d \times \mathcal{X} \times \Omega \to \N \cup \{\infty\} \]
where $\Omega$ is a probability space.  In this case, we require $T(\cdot,\cdot,\cdot,\omega)$ to satisfy properties (a)-(d) for almost all $\omega \in \Omega$.

If $u,v$ are functions on $\Z^d$, write $u \leq v$ if $u(x) \leq v(x)$ for all $x\in \Z^d$.  The discrete Laplacian $\Delta u$ of $u$ is the function 
    \begin{equation} \label{thediscretelaplacian} \Delta u(x) = \sum_{y\sim x} u(y) - 2d\, u(x) \end{equation}
where the sum is taken over the $2d$ lattice neighbors of $x$.  

Given a toppling procedure~$T$ and initial configuration $\eta$, the resulting configuration at time~$t$ is
	\[ \eta_t := \eta + \Delta T(t,\cdot,\eta). \]
We say that~$T$ is \emph{legal} for $\eta$ if for all $x\in \Z^d$ and all~$t\geq 0$ such that~$x$ topples at time~$t$, we have
	\[ \lim_{s \uparrow t} \eta_s(x) \geq 2d. \]
That is, in a legal toppling procedure only unstable sites are toppled.

We set
	\[ T(\infty,x,\eta) = \sup_{t\geq 0} T(t,x,\eta) \in \N \cup \{\infty\}. \]
We say that $T$ is \emph{finite} for initial configuration $\eta$, if $T(\infty,x,\eta)<\infty$ for all $x\in \Z^d$.  In this case, we say that $T$ is \emph{stabilizing} for $\eta$ if
	\[ \eta_{\infty} := \eta + \Delta T(\infty,\cdot,\eta) \leq \overline{2d-1}. \]
That is, in the final configuration $\eta_{\infty}$ no site is unstable.

The following lemma has been proved a number of times in various settings \cite{BLS,dhar,DF, FMR, MRZ}.

\begin{lemma} (Abelian property) 
For any $\eta \in \mathcal{X}$, if $T$ is a finite legal stabilizing toppling procedure for $\eta$, then any legal toppling procedure is finite for $\eta$.  Moreover if $T'$ is another legal stabilizing toppling procedure for $\eta$, then for all $x \in \Z^d$
	 \[ T(\infty,x,\eta) = T'(\infty,x,\eta). \]
\end{lemma}

If~$\eta$ is a particle configuration for which there exists a finite legal stabilizing toppling procedure, then we say that $\eta$ \emph{stabilizes}; otherwise, we say that~$\eta$ is \emph{exploding}.  If~$\eta$ stabilizes, then the function $u: \Z^d \to \N$ given by
	\begin{align} \label{theodometer}  u(x) &= T(\infty,x,\eta),
%	\\ &= \mbox{number of times $x$ topples during the stabilization of $\eta$}
 	\end{align}
where $T$ is any legal stabilizing toppling procedure for $\eta$, is called the \emph{odometer} of $\eta$. 
%The term odometer will always refer to the number of topplings occuring in a \emph{legal} toppling sequence.

Note that if every site topples at least once, then~$\eta$ must be exploding.  Otherwise, by the no infinite backward chain condition (d), some site~$x \in \Z^d$ must finish toppling no later than all of its neighbors do; since each neighbor topples at least once more,~$x$ receives~$2d$ additional particles and must topple again, a contradiction.  Thus we have shown

\begin{lemma} \cite[Theorem 2.8, item 4]{FMR}
\label{easyexplosion}
%If $u(x) \geq 1$ for all $x \in \Z^d$, then $\eta$ is exploding.
If $\eta$ stabilizes, then $u(x)=0$ for some $x \in \Z^d$.
\end{lemma}

%\begin{proof}
%Fix a legal stabilizing toppling procedure $T$.
%For each $x \in \Z^d$, let $t_x$ be the last time $x$ topples (or $t_x=0$ if $x$ never topples).  
%By the ``no infinite backward chain'' condition (d), there exists a site $x \in \Z^d$ with $t_x \leq t_y$ for all neighbors $y\sim x$.  Since $T$ is stabilizing, $x$ can receive at most $2d-1$ chips after time $t_x$, so some neighbor $y$ must never topple.  
%\end{proof}

%Thus in order to show that a configuration $\eta$ is exploding, it suffices to find a legal toppling procedure in which every site $x \in \Z^d$ topples at least once.  We will use this fact in section~\ref{sec:robustandexplosive}.

Let $\eta$ be a particle configuration that stabilizes, and let $u$ be its odometer function (\ref{theodometer}).   In our application, we will take $\eta = \sigma + n\delta_o$, where $\sigma$ is a robust background.  Since $\Delta u(x)$ counts the net number of particles exiting the site~$x$, the stabilization $\eta_{\infty}$ of $\eta$ is given by
	\[ \eta_{\infty} = \eta + \Delta u. \]  

\begin{definition}
Given a particle configuration $\eta$ on $\Z^d$, a function $u_1 : \Z^d \to \Z$, is called \emph{stabilizing} for~$\eta$ if 
	\[ \eta + \Delta u_1 \leq \overline{2d-1}. \]
\end{definition}
	
Informally, we may think of $\eta + \Delta u_1$ as the configuration obtained from $\eta$ by performing $u_1(x)$ topplings at each site $x \in \Z^d$.  Note, however, that the above definition makes no requirement that these topplings be legal; that is, they may produce sites with a negative number of particles.
	
Our proof of Theorem~\ref{thecube} rests on the following lemma, which characterizes the odometer function $u$ as minimal among all nonnegative stabilizing functions.  Deepak Dhar has aptly called this a ``least action principle,'' in the sense that the number of topplings in a legal toppling sequence is the minimum number required to stabilize the configuration.  In fact, more is true: not only is the total number of topplings minimized, but \emph{each vertex} does the minimum amount of work required of it to stabilize the configuration.  

According to the abelian property, if we use a legal toppling procedure to stabilize $\eta$, then each site $x$ topples exactly $u(x)$ times, regardless of the choice of procedure.  The least action principle says that in any sequence of topplings that stabilizes $\eta$, \emph{even if some of those topplings are illegal}, each site $x$ topples at least $u(x)$ times.
%as many times as in a legal stabilizing toppling sequence.

\begin{lemma} (Least Action Principle)
\label{odometerisminimal}
Let $\eta$ be a particle configuration on $\Z^d$ that is not exploding, and let $u$ be its odometer.
If $u_1 : \Z^d \to \N$ is stabilizing for $\eta$, then $u_1 \geq u$.
\end{lemma}

\begin{proof}
To compare $u_1$ to the odometer, we use the following discrete time legal toppling procedure $T'$.
Enumerate the sites in $\Z^d$.
Call a site~$x \in \Z^d$ \emph{ready} if it has at least~$2d$ particles and has toppled fewer than~$u_1(x)$ times.  At each time step, if there are any ready sites, topple the smallest ready site.
%Perform legal topplings sequentially, in any order,
%%(satisfying the requirements of Definition 2.1 of \cite{FMR}) 
%without allowing any site~$x$ to topple more than $u_1(x)$ times, until every site~$x$ either is stable, or has toppled exactly $u_1(x)$ times.

Write $u'(x)=T'(\infty,x,\eta)$ for the number of times $x$ topples during this procedure.  We will show that $u'=u$.  If $\eta' = \eta + \Delta u'$ is stable, then $T'$ is stabilizing as well as legal, so $u' = u$ by the abelian property.  Otherwise, $\eta'$ has some unstable site $y$.  We must have $u'(y) = u_1(y)$; otherwise, $y$ would still be ready.
Writing $u'' = u_1 - u'$, we obtain
	\[ (\eta + \Delta u_1)(y) = \eta'(y) + \Delta u''(y) \geq \eta'(y) \geq 2d \]
since $u''(y)=0$. This contradicts the assumption that $u_1$ is stabilizing.
\end{proof}
%\begin{proof}
%Let $x_1, \ldots, x_k$ be a legal toppling sequence that stabilizes $\sigma$.  Suppose for a contradiction that $u \not\leq u_1$, and choose $j$ minimal such that
%	\[ \# \{i<j \,:\, x_i = x_j\} = u_1(x_j).  \]
%By the minimality of $j$, each neighbor $y\sim x_j$ appears at most $u_1(y)$ times among $x_1, \ldots, x_{j-1}$, so the number of particles present at $x_j$ after toppling $x_1, \ldots, x_{j-1}$ is at most
%	\[ \sigma(x_j) - 2d u_1(x_j) + \sum_{y \sim x_j} u_1(y) \leq 2d-1. \]
%This contradicts the legality of toppling $x_j$. 
%\end{proof}

We pause here to record a closely related fact.  
If $u_1, u_2$ are functions on $\Z^d$, write $\min(u_1,u_2)$ for their pointwise minimum.
If $x\in \Z^d$ is a site where $u_1(x) \leq u_2(x)$, then 
\begin{eqnarray*}
\Delta\min(u_1,u_2)(x) & = & \sum_{y\sim x} \min(u_1(y),u_2(y)) - 2d\,u_1(x)\\
 & \leq &  \sum_{y\sim x} u_1(y) - 2d\,u_1(x) = \Delta u_1(x).
\end{eqnarray*}
Likewise, if $u_1(x) > u_2(x)$, then $\Delta\min(u_1,u_2)(x) \leq \Delta u_2(x)$. So 
	\[ \Delta\min(u_1,u_2) \leq \max(\Delta u_1,\Delta u_2). \]
As a consequence, we obtain the following.

\begin{lemma} 
If $u_1$ and $u_2$ are stabilizing for $\sigma$, then $\min(u_1, u_2)$ is also stabilizing for $\sigma$.
\label{tropical}
\end{lemma}

\begin{proof}
%Since $\Delta \min(u_1,u_2) \leq \max (\Delta u_1, \Delta u_2)$, we have
$\sigma + \Delta \min(u_1,u_2) \leq \max(\sigma+\Delta u_1, \sigma+ \Delta u_2) \leq \overline{2d-1}. $
\end{proof}

The set of stabilizing functions is also closed under adding any constant function, giving it the structure of a module over the tropical semiring $(\Z,\min,+)$.  A related module is studied in~\cite{HMY}.

\section{Growth Rates}
\label{growthratessection}

Fix an integer $h \leq 2d-2$, and let $\eta$ be the configuration $\overline{h}+n\delta_o$ on~$\Z^d$.  Let $S_n$ be the set of sites that ever topple or receive a particle during the stabilization of~$\eta$ (in~\cite{LP09a} these were called ``visited'' sites).  Note that if~$y$ receives a particle, then one of its neighbors must have toppled.  Thus $S_n$ is related to the set $T_n$ of sites that topple by
%	\[ T_n  \subset S_n \subset T_n \cup \partial T_n  \]
	\[ S_n = T_n \cup \partial T_n \]
where for $A \subset \Z^d$ we write 
	\[ \partial A = \{y \in \Z^d \,:\, y\not \in A, \exists z\in A, z\sim y\}. \]

Write $|x|=(x_1^2+\ldots+x_d^2)^{1/2}$ for the Euclidean norm on $\Z^d$, and for~$r>0$ let 
	\[ B_r = \{x \in \Z^d \,:\, |x|<r\} \]
be the ball of radius~$r$ centered at the origin in~$\Z^d$.  Let $\omega_d$ be the volume of the unit ball in~$\R^d$.
For the proof of Theorem~\ref{thecube}, we take as a starting point the following result of \cite{LP09a}.

\begin{theorem} \cite[Theorem~4.1]{LP09a}
\label{lowheightbounds}
Fix an integer $h \leq 2d-2$.  For any $n\geq 1$, we have
	\[ B_{c_1r - c_2} \subset S_n \]
where $r$ is such that $n = \omega_d r^d$, and
	\[ c_1 = (2d-1-h)^{-1/d} \]
and $c_2$ is a constant depending only on $d$.  Moreover if $h \leq d-1$, then for any $n \geq 1$ and any $\epsilon>0$ we have
	 \begin{equation} \label{lowheightouterbound} S_n \subset B_{c'_1r + c'_2} \end{equation}
where
 	\[ c'_1 = (d-\epsilon-h)^{-1/d} \]
and $c'_2$ is independent of $n$ but may depend on $d$, $h$ and $\epsilon$.
\end{theorem}

Note that $h$ may be negative, in which case the background $\overline{h}$ corresponds to each site in $\Z^d$ starting with a ``hole'' of depth $H=-h$. 

We are grateful to Haiyan Liu for pointing out a gap in the proof of the outer bound (\ref{lowheightouterbound}).  The gap occurs in Lemma~4.2 of~\cite{LP09a}, which is valid only for $H \geq 0$.  We correct this gap in section~\ref{sec:correction}.  Next, in section~\ref{sec:mainproof}, we explain our technique of ``background modification,'' and use it to deduce Theorem~\ref{thecube} from Theorem~\ref{lowheightbounds}.

\subsection{Low Background Height}
\label{sec:correction}
Fix $0 \leq h \leq d-1$, let $\eta$ be the configuration $\overline{h}+n\delta_o$ on $\Z^d$, and consider the odometer function
	\[ u_n(x) = \mbox{number of times $x$ topples during the stabilization of } \eta. \]
The normalization of the odometer function and of the discrete Laplacian (\ref{thediscretelaplacian}) differs by a factor of $2d$ from the one used in~\cite{LP09a}.  It is the most convenient normalization for the abelian sandpile, since $2d$ particles move in every toppling.

In \cite{LP09a} it is proved that for every site $x\in \Z^d$ with $c'_1 r -1 < |x| \leq c'_1 r$ we have
	\[ u_n(x) \leq c \]
where $c$ is a constant which may depend on $d$, $h$ and $\epsilon$ but not on $n$.  (In the notation of \cite{LP09a}, $c=c'_2/2d$.)

It follows that $u_n$ is uniformly bounded outside the ball $B_{c'_1 r}$; indeed, if $|x|>c'_1 r$, then setting $n' = \ceiling{\omega_d (|x|/c'_1)^d}$, since $n \leq n'$ we have by the abelian property
	\[ u_n(x) \leq u_{n'}(x) \leq c.  \]
The next lemma shows that in fact, $u_n=0$ outside the slightly larger ball $B_{c'_1 r + c-1}$.  Hence $T_n \subset B_{c'_1 r + c -1}$, and hence $S_n \subset B_{c'_1r + c}$, which completes the proof of (\ref{lowheightouterbound}).

\begin{lemma}
For all $j=0,1,\ldots,c$ and all $x \in \Z^d$ with $|x|>c'_1r + j-1$, we have
	\[ u_n(x) \leq c-j. \]
\end{lemma}

\begin{proof}
Let
	\[ R_j = \{x \in \Z^d \,:\, |x| > c'_1r + j -1 \}. \]
Note that for any $x \in R_j$, all neighbors $y\sim x$ lie in $R_{j-1}$, and at least~$d$ neighbors have $|y| \geq |x|$, so at least~$d$ neighbors lie in $R_j$.  

We will prove the lemma by induction on $j$.  Let
	\[ U_j = \{x \in R_j \,|\, u_n(x) > c-j \}. \]
If $U_j$ is empty, then the proof is complete.  Otherwise, by the no infinite backward chain condition, there exists a site $x\in U_j$ that finishes toppling no later than all of its neighbors in~$U_j$.  By the inductive hypothesis, every neighbor $y$ of $x$ satisfies.
	\[ u_n(y) \leq c-j+1. \]
Just before $x$ topples for the last time, each neighbor~$y \in U_j$ has not yet toppled for the last time, so $y$ has toppled at most~$c-j$ times.  Moreover, each neighbor~$y \in R_j - U_j$ has toppled at most~$c-j$ times; and each neighbor $y \notin R_j$ has toppled at most~$c-j+1$ times.  Hence, just before it topples for the last time, $x$ has received at most~$d(c-j) + d(c-j+1)$ chips and emitted at least~$2d(c-j)$ chips, leaving it with at most~$h+d$ chips.  
%So just before $x$ topples for the last time, it has
%	\[ h - 2d(c-j) + d(c-j+1) + d(c-j) = h+d \]
%chips.  
Since~$h\leq d-1$, this is not enough chips to topple, which gives the required contradiction.
\end{proof}

\subsection{High Background Height}
\label{sec:mainproof}
To prove Theorem \ref{thecube} using the least action principle (Lemma \ref{odometerisminimal}), for each coordinate $i=1,\ldots,d$ we will construct a toppling function~$g_i$ supported in the slab
    \begin{equation} \label{theslab} A_{i,r} = \{x \in \Z^d \,:\, |x_i| \leq r \}. \end{equation}
%    for the infinite slab of width $2r+1$ centered at the origin in $\Z^d$.
The effect of toppling according to~$g_i$ will be to modify the constant background height $h$ by ``clearing out'' particles down to height at most $d-1$ in a smaller slab $A_{i,r_0}$ and ``piling them up'' to height at most $2d-1$ outside $A_{i,r_0}$.
We will see that this can be done while keeping $r_0$ proportional to $r$.
 
On this modified background, $n$ particles at the origin will spread with a growth rate at most according to $h = d-1$, provided $n$ is small enough so that the particles do not spread outside $Q_{r_0}$.  This growth rate is controlled by Theorem~\ref{lowheightbounds}: $n$ particles on constant background height $d-1$ in $\Z^d$ spread at most a distance of order $n^{1/d}$.  Since $r_0$ is proportional to $r$, we can therefore choose $n$ proportional to $r^d$.
 
The desired background modification can be accomplished a function of just one coordinate, $g_i(x_1,\ldots,x_d) = g(x_i)$.  The next lemma spares us the need to specify $g$ explicitly; it suffices to specify how the background is modified.  In the lemma, $g$ plays the role of toppling function on $\Z$, and $f$ represents the net change in height of the configuration. The conditions (\ref{conservationandcenterofmass}) mean in words that topplings cannot change the total number of particles, nor the center of mass of a configuration.

\begin{lemma}
\label{onedimensional}
If $f: \Z \to \Z$ is supported on a finite interval $I = [-a,b]$, and 
	\begin{equation} \label{conservationandcenterofmass} \sum_{y \in I} f(y) = \sum_{y \in I} y f(y) = 0, \end{equation}
then $f=\Delta g$ for an integer-valued function $g$ supported on the interval $I'=[1-a,b-1]$.  Moreover, if there are no $x_1<x_2<x_3$ such that $f(x_1)<0$, $f(x_2)>0$ and $f(x_3)<0$, then $g\geq 0$.
%Moreover, if the set $J = \{x \in \Z \,:\, f(x)<0 \}$ is an interval, then $g \geq 0$.
%this condition is easier to parse, but slightly more restrictive.
\end{lemma}

\begin{proof}
Let 
	\[ g(x) = \sum_{y=-a}^{x-1} (x-y)f(y). \]
Then for $x \geq b$ we have
	\[ g(x) = x \sum_{y=-a}^b f(y) - \sum_{y=-a}^b y f(y) = 0 \]
so $g$ is supported on $I'$.  Also
	\begin{align*} \Delta g(x) &= g(x+1) - 2g(x) + g(x-1) \\
					&= f(x) + \sum_{y=-a}^{x-2} ((x+1-y) - 2(x-y) + (x-1-y)) f(y) \\
					&= f(x)
					\end{align*}
as desired.  

If $g(z)<0$ for some $z$, then since $g(-a)=g(b)=0$, the difference
	\[ Dg(y) = g(y+1) - g(y) \]
satisfies $Dg(y_1)<0$ and $Dg(y_2)>0$ for some $y_1< z \leq y_2$.  Hence the second difference
	\[ f(x) = \Delta g(x) = Dg(x) - Dg(x-1) \]
%	\[ D^2g(y) = Dg(y+1)-Dg(y) \]
satisfies $f(x_1)<0$, $f(x_2)>0$ and $f(x_3)<0$ for some $x_1 \leq y_1 < x_2 \leq y_2 < x_3$.  
\end{proof}

%Write
%	\[ Q_r = \{x\in \Z^d \,:\, \max |x_i| \leq r \} \]
%for the cube of side length $2r+1$ centered at the origin in $\Z^d$. 

\begin{proof}[Proof of Theorem~\ref{thecube}]
For each $i=1,\ldots,d$ we will construct a nonnegative function $u_i$ on $\Z^d$ which is stabilizing for the configuration $\overline{h} + n \delta_o$, and supported on the infinite slab $A_{i,r}$; see (\ref{theslab}).

By the least action principle, Lemma~\ref{odometerisminimal}, the odometer function $u$ satisfies $u \leq u_i$ for $i=1,\ldots, d$.  Since $T_{n,d,\overline{h}}$ is the support of $u$, we obtain
    \[ T_{n,d,\overline{h}} \subseteq \bigcap_{i=1}^d A_{i,r} = Q_r. \]
%since each of the functions $u_i$ dominates the odometer.

To construct $u_i$, let $w : \Z^d \to \N$ be the odometer function for the configuration $\overline{d-1}+n\delta_o$.  By Theorem~\ref{lowheightbounds}, if $n$ is sufficiently large, then $w$ is supported on the ball centered at the origin of radius
    \begin{equation} \label{rho} \rho =  \left(1+\frac{\epsilon}{2d}\right) \left( \frac{n}{\omega_d} \right)^{1/d}. 
\end{equation}
In particular, $w$ vanishes outside the cube $Q_{\rho}$.

Let $r_0$ be the smallest integer multiple of $2d-1-h$ exceeding $\rho$, and let
    \[ r_1 = \frac{d}{2d-1-h} r_0. \]
Let $f: \Z \to \Z$ be given by
    \[ f(x) = \begin{cases} 2(d-1-h), & x=0 \\
                     d-1-h, & 0 < |x| < r_0 \\
                     2d-1-h, & r_0 \leq |x| < r_1 \\
                          0, & |x|\geq r_1.
                \end{cases} \]
Then with $I = [1-r_1,r_1-1]$
    \begin{align*} \sum_{y \in I} f(y) &= 2r_0(d-1-h) + (2r_1-2r_0)(2d-1-h) \\
                            &= -2dr_0 + 2(2d-1-h)r_1 = 0. \end{align*}
Since $f(y)=f(-y)$ we have $\sum_{y \in I} yf(y) =0$.  By Lemma~\ref{onedimensional}, $f = \Delta g$ for a nonnegative integer-valued function $g$ supported on the interval $I'=[2-r_1,r_1-2]$.

For $x = (x_1, \ldots, x_d) \in \Z^d$, define
    \[ u_i(x) = w(x) + g(x_i). \]
Note that the function $g(x_i)$ has $d$-dimensional Laplacian $f(x_i)$.
Inside the cube $Q_\rho$, we have $f(x_i) \leq d-1-h$, hence inside $Q_\rho$
    \[ h + n\delta_o + \Delta u_i \leq d-1 + n\delta_o + \Delta w \leq 2d-1. \]
Outside $Q_\rho$, since $w$ vanishes and $f(x_i) \leq 2d-1-h$, we have
    \[ h + n\delta_o + \Delta u_i \leq 2d-1. \]
Thus $u_i$ is stabilizing for $h+ n\delta_o$.  Moreover, since
    \begin{align*} r_1 &\leq d \left( \frac{\rho}{2d-1-h} + 1 \right) \\
                    &= \frac{d + \epsilon/2}{2d-1-h} \left( \frac{n}{\omega_d} \right)^{1/d} + d \end{align*}
we have $r_1 \leq r$ for sufficiently large $n$, hence $u_i$ is supported on the slab $A_{i,r}$ as desired.
\end{proof}

\begin{figure}
\centering
\includegraphics[width=.48\textwidth]{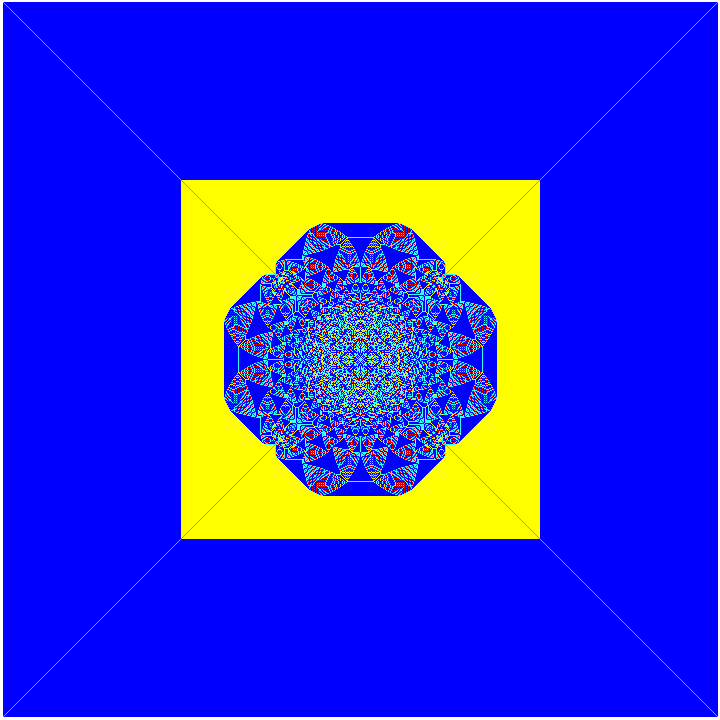}
\includegraphics[width=.48\textwidth]{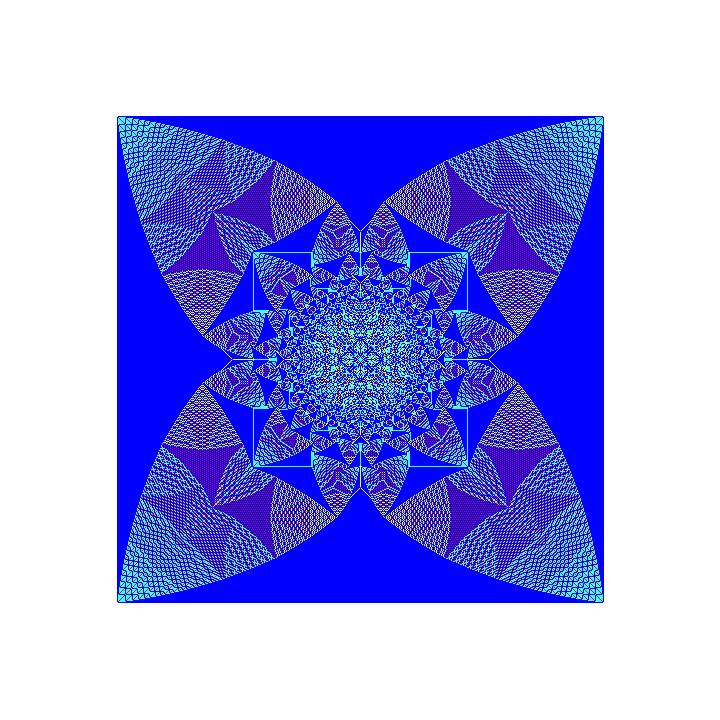}
\caption{Left: The stable configuration $\overline{2}+n\delta_o+\Delta \min(u_1,u_2)$ constructed in the proof of Theorem~\ref{thecube}.  Sites with negative height, along the diagonals of the square, are colored orange. Right: The stabilization of $\overline{2}+n\delta_o$.  Here $n=10^5$.  
}
\label{stabilizingodom}
\end{figure}

We remark that in addition to bounding the set of sites $T_{n,d,\overline{h}}$ that topple, the proof gives a bound on the odometer function
    \[ u(x) 
    %u_{n,d,h}(x)   (subscript notation clashes with u_i above and below)
    = \#\, \mbox{times $x$ topples in the stabilization of $\overline{h}+n\delta_o$ in $\Z^d$}, \]
namely
    \begin{align*} u(x) &\leq \min (u_1(x),\ldots,u_d(x)) \\
        &= w(x) + \min (g(x_1), \ldots, g(x_d)) \\
        &= w(x) + g(\max |x_i|). \end{align*}
By Lemma~\ref{tropical}, the right side is stabilizing for $\overline{h}+n\delta_o$.  The resulting stable configuration in the case $d=h=2$ is pictured in Figure~\ref{stabilizingodom}.

The proof of Theorem~\ref{rowsandcolumns} is identical to that of Theorem \ref{thecube}, except that we now choose
\[ f(x) = \begin{cases} 2-2d, & x=0 \\
                  -d, & 0 < |x| < r_0 \mbox{ and } m\nmid x \\
                1-d, & 0 < |x| < r_0 \mbox{ and } m|x \\
                 0, & r_0 \leq |x| < r_1 \mbox{ and } m\nmid x \\
                 1, & r_0 \leq |x| < r_1 \mbox{ and } m|x \\
                 0, & |x|\geq r_1,
                \end{cases} \]
where $r_0$ the smallest integer exceeding $\rho =  (1+\frac{\epsilon}{2d}) \big( \frac{n}{\omega_d} \big)^{1/d}$,
%(see equation \ref{rho})
and 
\[ r_1 = m(dr_0 - 1). \]

Then we have that 
\[
r_1 \leq m \left(d + \frac{\epsilon}{2}\right) \left( \frac{n}{\omega_d} \right)^{1/d} + m(d-1), 
\]
so that again $r_1 \leq r$ for sufficiently large $n$.

\section{Robust and Explosive Backgrounds}
\label{sec:robustandexplosive}

Write $\psi_1 = e_1, \ldots, \psi_d = e_d, \psi_{d+1}=-e_1, \ldots, \psi_{2d} = -e_d$
for the $2d$ coordinate directions in $\Z^d$.  
If $R$ is a rectangular prism in $\Z^d$, write
	\[ F_i(R) = \{y: y \not \in R, y-\psi_i \in R\} \]
for the outer face of $R$ in direction $\psi_i$.

We will deduce Propositions~\ref{latticebootstrapintro} and~\ref{randombootstrapintro} from the following slightly more general result. 
%The method of proof closely resembles Straley's argument for bootstrap percolation in \cite{enter}.

\begin{theorem}
\label{disaster}
Let $\sigma$ be a background on $\Z^d$ satisfying
	\begin{itemize}
	\item[{\em (i)}] $\sigma(x) \geq 2d-2$ for all $x \in \Z^d$; and
	\item[{\em (ii)}] There exists $r_0 \in \N$ such that for all $r \geq r_0$, each outer face $F_i(Q_r)$ contains a site $x$ with $\sigma(x) \geq 2d-1$.
	\end{itemize}
Then $\sigma$ is explosive.
\end{theorem}

\begin{proof}
By Lemma~\ref{easyexplosion}, in order to prove that a configuration~$\eta$ on~$\Z^d$ is exploding, it suffices to find a toppling procedure in which every site in~$\Z^d$ topples at least once.

From \cite[Theorem~4.1]{shapes}, if the background height is exactly $2d-2$, then for every $n$, the set of sites that topple during stabilization forms a cube $Q_r$, and we can
%and every boundary site of this cube topples exactly once. Therefore, every site on the outer boundary of the cube has $2d-1$ particles.  
choose $n$ so that $r \geq r_0$.

%We add $n$ particles at the origin, and organize the topplings in stages. 
%In the first stage, we perform only the topplings that would occur if the background height were $2d-2$. Then we examine the outer boundary of the cube of toppled sites. If a site on the outer boundary in fact had height $2d-1$, then it is now unstable.

Let $R_0 = Q_r$ and 
	\[ R_k = R_{k-1} \cup F_{k~\mbox{\scriptsize mod } 2d}(R_{k-1}), \qquad k\geq 1. \]
We will define a toppling order in stages $k=0,1,2,\ldots$ so that at the end of stage $k$, all sites in $R_k$ have toppled at least once, and no other sites have toppled.  Since $\bigcup_{k\geq 0} R_k = \Z^d$, it follows that every site in~$\Z^d$ topples at least once, so $\sigma + n \delta_o$ is exploding.

During stage $0$, we perform all the topplings that occur in the stabilization of $\overline{2d-2}+n\delta_o$. Then every site in the cube $R_0$ has toppled at least once, and no other sites have toppled.  Hence by (i), every site in every outer face of $R_0$ now has at least $2d-1$ particles, and by (ii), in every outer face there is at least one unstable site. 

%In stage $1$, we choose a direction $i$, and topple every site in $F_{i,Q_r}$. We can do this by starting a wave at the unstable site. Now we have that all sites in $R_1 = Q_r \cup F_{i,Q_r$ have toppled at least once, and no site outside $R_1$ has toppled.

The remaining stages are defined inductively.
After stage $k-1$, every site in $F = F_{k~\mbox{\scriptsize mod } 2d}(R_{k-1})$ has at least $2d-1$ particles, and at least one site $y \in F$ is unstable.  Topple first $y$, then its neighbors in $F$, then the sites in $F$ at distance $2$ from $y$, and so on, until all sites in $F$ have toppled once.  Now every site in $R_k$ has toppled at least once, and no sites outside $R_k$ have toppled, completing the inductive step.  
\end{proof}

To deduce Proposition~\ref{latticebootstrapintro} from Theorem~\ref{disaster}, since 
	\[ \gcd(x_{1j},\ldots,x_{dj})=1, \qquad j=1,\ldots,d \]
there exist integers $a_{ij}$ for $1 \leq i,j \leq d$, satisfying
	\[ \sum_{i=1}^d a_{ij} x_{ij} = 1, \qquad j=1,\ldots,d.  \]
Then for each $j=1,\ldots,d$, the vector
	\[ \mathbf{v}_j = \sum_{i=1}^d a_{ij} \xx_i \in L \]
has $e_j$-coordinate $v_{jj}=1$, so any hyperplane in $\Z^d$ parallel to one of the coordinate hyperplanes intersects $L$.  Moreover, $L$ contains the vectors~$De_j$ for $j=1,\ldots,d$, where $D=\det(x_{ij})_{i,j=1}^d \neq 0$.  Thus $L$ intersects every face $F_i(Q_r)$ when $r \geq |D|/2$.

To deduce Proposition~\ref{randombootstrapintro}, it remains to check that the configuration $\overline{2d-2}+\beta$ on $\Z^d$ satisfies condition (ii) of Theorem~\ref{disaster} with probability~$1$.  Write $\mathcal{E}_{i,r}$ for the event that $\beta(x)=0$ for all $x \in F_i(Q_r)$.  By the independence of the Bernoulli random variables $\beta(x)$, this event has probability
	\[ \PP(\mathcal{E}_{i,r}) = (1- \epsilon)^{|F_i(Q_r)|} \leq (1-\epsilon)^r. \]
In particular, $\sum_{r \geq 1} \sum_{i=1}^{2d} \PP(\mathcal{E}_{i,r}) < \infty$.  By the Borel-Cantelli lemma, with probability $1$ only finitely many of the events $\mathcal{E}_{i,r}$ occur.   We remark on the similarity between this argument and Straley's argument for bootstrap percolation \cite{enter}.  

%For the configuration $2+1_L+n\delta_o$, we define a toppling order in stages so that at the end of stage $k$, all sites in $Q_{r+k}$ have toppled at least once, and no other sites have toppled.  During stage $0$, perform all the topplings that occur in the stabilization of $2+n\delta_0$.  The remaining stages are defined inductively.  Note that after stage $k$, all sites in $B_{r+k}$ except for the four corners have height at least $3$; the corners have height at least $2$; and each of the four sides has at least one site belonging to $L$, which has height at least $4$.  During stage $k+1$, topple every site in $B_{r+k}$ exactly once, beginning with the sites of height $4$ and ending with the corners.  After these topplings, all sites in $Q_{r+k+1}$ have toppled at least once, and no other sites have toppled, so the inductive step is complete.

%Since every site in $\Z^2$ topples, we conclude that $2+1_L+n\delta_0$ is exploding.
%\end{proof}

\medskip

We define the box $\mathcal{B}_r$ as 
\[ \mathcal{B}_r =  \partial Q_r = \bigcup_{i=1}^{2d} F_i(Q_r). \]
The following theorem is a partial converse to Theorem~\ref{disaster}, and gives a counterexample to Remark 7.1 in \cite{feyredig}. 

\begin{theorem}
\label{thm:boxes}
%Let $\mathcal{B}_1, \mathcal{B}_2, \ldots$ be a sequence of disjoint nested boxes in $\Z^d$.
Let $r_1, r_2, \ldots$ be an increasing sequence of positive integers.
Let $\sigma$ be a stable background on $\Z^d$ in which every site in $\mathcal{B}_{r_1} \cup \mathcal{B}_{r_2} \cup \ldots$ has at most $2d-2$ particles.

Then $\sigma$ is robust.
\end{theorem}

%Note that the background in Theorem~\ref{rowsandcolumns} satisfies the hypotheses of this theorem.
The proof uses the following lemma, which follows from \cite[Lemma~4.2]{shapes} and the abelian property.

\begin{lemma}
\label{onebox}
If $\sigma$ is a stable background satisfying $\sigma(x) \leq 2d-2$ for all $x \in \mathcal{B}_r$, then no sites outside $Q_r$ topple during the stabilization of $\sigma + \delta_o$.
\end{lemma}

\begin{proof}[Proof of Theorem~\ref{thm:boxes}]
We need to show that $\sigma + n\delta_o$ stabilizes in finitely many topplings, for every $n \in \N$.
We induct on $n$ to show that no sites outside $Q_{r_n}$ topple during the stabilization of $\sigma + n \delta_o$.

By Lemma~\ref{onebox}, no sites outside $Q_{r_1}$ topple during the stabilization of $\sigma + \delta_o$. 

Let $\sigma_n$ be the stabilization of $\sigma + n\delta_o$.  By the inductive hypothesis, no sites topple outside $Q_{r_n}$ during this stabilization, so $\sigma_n(x) \leq 2d-2$ for all $x \in \mathcal{B}_{r_{n+1}}$.  By Lemma~\ref{onebox}, no sites outside $Q_{r_{n+1}}$ topple during the stabilization of $\sigma_n + \delta_o$.  By the abelian property, a site topples during the stabilization of $\sigma + (n+1) \delta_o$ if and only if it topples during the stabilization of $\sigma + n\delta_o$ or during the stabilization of $\sigma_n + \delta_o$.  This completes the inductive step.
\end{proof}

\begin{remark}
Theorem~\ref{thm:boxes} remains true for arbitrary disjoint rectangular boxes surrounding the origin; they need not be cubical or centered at the origin. 
\end{remark}

\section{Dimensional Reduction}

\begin{figure}
\centering
\includegraphics[width=.48\textwidth]{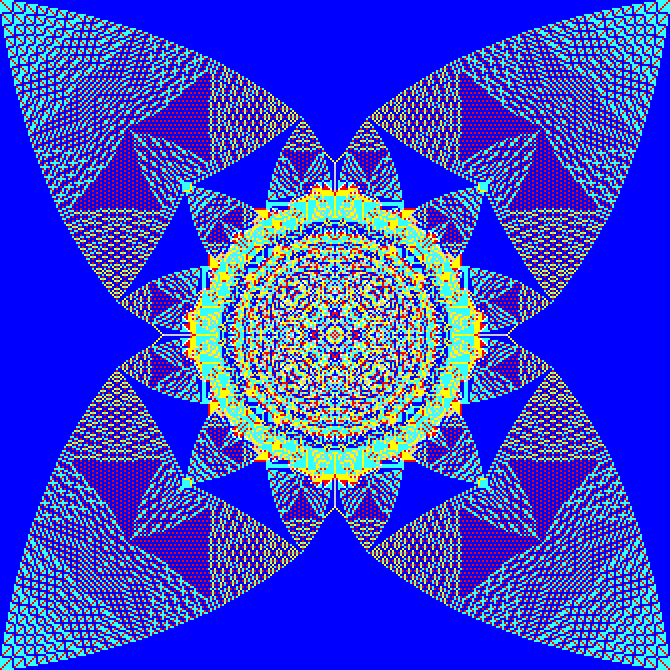}
\includegraphics[width=.48\textwidth]{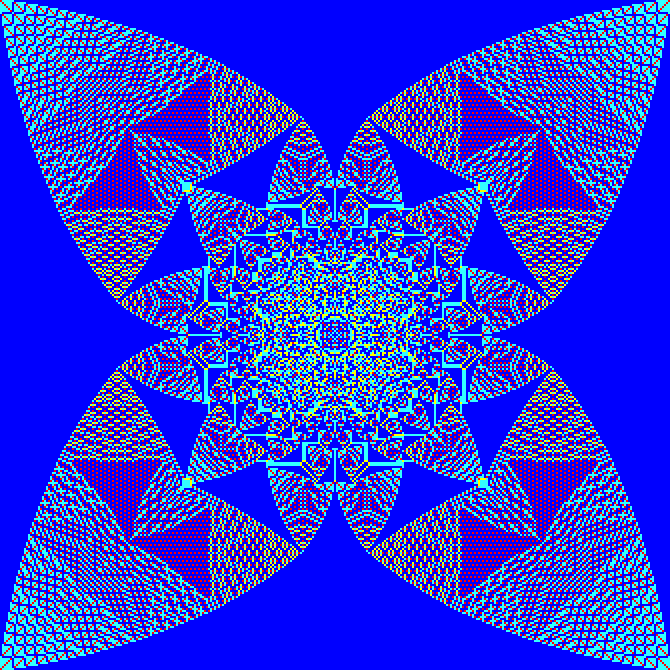}
%the 2-dim config is "full" - adding one more particle will cause the square to expand.
\caption{Left: A two-dimensional slice through the origin of the sandpile of $n=5\cdot 10^6$ particles in $\Z^3$ on background height $h=4$.  Right: The sandpile of $m=47465$ particles in $\Z^2$ on background height $h=2$.  Color scheme on left: sites colored blue have $5$ particles, turquoise $4$, yellow $3$, red $2$, gray $1$, white $0$.  On right: blue $3$ particles, turquoise $2$, yellow $1$, red $0$.}
\label{fig:dimensionalreduction}
\end{figure}

%\begin{figure}
%\centering
%\includegraphics[width=.241\textwidth]{sandpile3d5mslice2.png}
%\includegraphics[width=.241\textwidth]{sandpile3d5mslice5.png}
%\includegraphics[width=.241\textwidth]{sandpile3d5mslice8.png}
%\includegraphics[width=.241\textwidth]{sandpile3d5mslice10.png}
%\caption{A sequence of two-dimensional slices of the sandpile of $5\cdot 10^6$ particles in $\Z^3$ at $h=4$, from left to right in order of increasing distance from the origin.}
%\label{fig:slices}
%\end{figure}

Our argument used properties of the one-dimensional sandpile to bound the growth rate of higher-dimensional sandpiles.  There appears to be a deeper relationship between sandpiles in $d$ and $d-1$ dimensions, 
which we formulate in the following \emph{dimensional reduction conjecture}.  
For $x \in \Z^d$, let $\sigma_{n,d}(x)$ be the final number of particles present at $x$ in the stabilization of $\overline{2d-2} + n\delta_o$.   Write rad$(n,d)$ for the radius of the cube $\{x \in \Z^d| \sigma_{n,d}(x)>0\}$.  Note that by Theorem~\ref{thecube}, rad$(n,d)$ has order $n^{1/d}$.
%a literal reading of the thm gives only at most n^{1/d}

Identifying $\Z^{d-1}$ with the hyperplane $x_d=0$ in $\Z^d$, we would like to compare the slice through the origin  $Q_{\mbox{\scriptsize rad}(n,d)} \cap \Z^{d-1}$ of the $d$-dimensional sandpile started from $n$ particles with a $(d-1)$-dimensional sandpile started from some number $m$ of particles.  Given $m$, $n$ and $d$, let us call a site $x = (x_1,\ldots,x_{d-1}) \in \Z^{d-1}$ an \emph{exact match} if
	\[ \sigma_{n,d}(x_1,\ldots,x_{d-1},0) = 2+\sigma_{m,d-1}(x_1,\ldots,x_{d-1}). \]
Given $0<\lambda<1$, consider the subset $A_\lambda$ of the slice through the origin
	\[ A_\lambda = \left( Q_{\mbox{\scriptsize rad}(n,d)} - Q_{\lambda \mbox{\scriptsize rad}(n,d)} \right) \cap \Z^{d-1}. \]
	
\begin{conjecture}
There exists a constant $\lambda = \lambda_d < 1$ such that for all $n \geq 1$ there exists $m \geq 1$ such that all but $O(\mbox{\em rad}(n,d)^{d-2})$ sites in $A_\lambda$ are exact matches.
\end{conjecture}

The case $d=3$ is 
illustrated in Figure~\ref{fig:dimensionalreduction}.  Amazingly, except in a region near the origin, the two pictures shown in the figure agree pixel for pixel.  For some rare values of $n$, certain small ``defects'' or ``filaments'' in the two pictures fail to match exactly, which is why we exclude up to $O(\mbox{rad}(n,d)^{d-2})$ sites.  For simplicity, we have restricted our formulation to the slice through the origin, but dimensional reduction seems to occur in all slices except for those close to the boundary of the cube.  The value of $m$ is the same for all of these slices.  We first learned of the dimensional reduction phenomenon from Deepak Dhar.

%\clearpage

\end{document}